\def \NN {\mathbb N}

\def \RR {\mathbb R}

\def \A  {{\mathcal A}}

\def \H  {{\mathcal H}}

\def \I  {{\mathcal I}}
\def \J  {{\mathcal J}}

\def \P  {{\mathcal P}}

\def \M  {{\mathcal M}}

\documentclass[12pt,reqno]{amsart}

\usepackage{amsmath,amssymb}

\usepackage{url} 

\numberwithin{equation}{section}

\begin{document}


\title[]{An extension of the Bourgain-Sarnak-Ziegler  \\ theorem with modular applications}
\author[]{M.Cafferata, A.Perelli \lowercase{and} A.Zaccagnini}
\maketitle

{\bf Abstract.} We first prove an extension of the Bourgain-Sarnak-Ziegler theorem, relaxing some conditions and giving quantitative estimates. Then we apply our extension to bound certain exponential sums, where the coefficients come from modular forms and the exponential involves polynomial sequences of any degree.



\vskip.5cm
\section{Introduction}

A well known theorem by Bourgain-Sarnak-Ziegler \cite{BSZ/2013} (BSZ theorem for short), see also K\'atai \cite{Kat/1986} for an earlier version, asserts that given a small parameter $\tau>0$ and two arithmetical functions $a(n)$ and $\phi(n)$, with $|a(n)|\leq 1$ multiplicative and $|\phi(n)|\leq 1$ satisfying
\[
\Big|\sum_{m\leq M} \phi(pm) \overline{\phi(qm)}\Big| \leq \tau M
\]
for all primes $p,q\leq e^{1/\tau}$, $p\neq q$, and $M$ sufficiently large, then for $N$ large enough one has
\[
\Big|\sum_{n\leq N} a(n)\phi(n)\Big| \leq 2\sqrt{\tau\log(1/\tau)} N.
\]
The BSZ theorem has many interesting applications, typically in the framework of Sarnak's M\"obius Randomness conjecture \cite{Sar/2011}, where $a(n)=\mu(n)$ while $\phi(n)$ ranges from classical exponential cases to several new examples coming from dynamical systems.

\smallskip
In this paper we first establish an extension of the BSZ theorem which, essentially, includes multiplicative functions $a(n)$ that are suitably bounded on average. Then we apply it to bound certain polynomial exponential sums with modular coefficients. As it will be clear in a moment, such an extended BSZ theorem may be applied to a variety of other cases.

\smallskip
Throughout the paper $p$ denotes a prime number, $|\A|$ denotes the cardinality of a set $\A\subset \NN$, $f\asymp g$ means $f \ll g \ll f$ and an empty product equals 1. We prove an extension of the BSZ theorem under the following conditions.

\smallskip
{\bf Assumptions.} Let $x$ be sufficiently large, $H=H(x)$ and $K=K(x)$ be parameters satisfying
\begin{equation}
\label{1-1}
\log^{\delta} x < H < K < x^{\delta}
\end{equation}
with some $0<\delta\leq 1/10$, say, and let 
\[
\P=\{z<p\leq w\} \quad \text{and} \quad P=\prod_{p\in\P}p.
\]
Suppose that $a(n)$ is a multiplicative arithmetical function satisfying $a(p) \ll 1$ and $\phi(n)$ is a bounded arithmetical function. Moreover, suppose that the following assumptions are satisfied whenever
\[
H^2/2 \leq z < w \leq 2K^2:
\]

(a) if $\overline{P}=1$ or $\overline{P}=P$ and $y\gg x/w$, then as $x\to\infty$ we have
\[
\sum_{\substack{n\leq y \\ (n,\overline{P})=1}} |a(n)|^2 \ll y \prod_{p|\overline{P}}\Big(1-\frac{1}{p}\Big),
\]

(b) if $w-z \asymp \sqrt{z}$ and $y \asymp x/z$, then as $x\to\infty$ we have
\[
\sum_{\substack{p,q\in\P \\ p\neq q}} \Big| \sum_{m\leq y} \phi(pm) \overline{\phi(qm)} \Big| \ll \tau \frac{z y}{\log^2z} \quad \text{with some $\tau=\tau(x)\leq 1$},
\]
where the constants in the $\ll$-symbols may depend at most on $a(n)$, $\phi(n)$ and $\delta$. \qed

\smallskip
Note that $\tau$ in (b) represents, essentially, the saving over the trivial bound. Finally, let
\begin{equation}
\label{1-2}
S(x) = \sum_{n\leq x} a(n) \phi(n).
\end{equation}
The extension of the BSZ theorem is as follows.

\medskip
{\bf Theorem 1.} {\sl Under the above assumptions, as $x\to\infty$ we have
\[
S(x) \ll x \Big(\frac{1}{\sqrt{H\log H}} +  \sqrt{\tau} + \frac{\log H}{\log K}\Big),
\]
where the constant in the $\ll$-symbol depends at most on $a(n)$, $\phi(n)$ and $\delta$.}

\medskip
We remark that the assumptions in Theorem 1 may be somewhat relaxed.

\medskip
Turning to the applications to exponential sums, let $e(\theta) = e^{2\pi i\theta}$ and
\[
S_a(x,\xi) = \sum_{n\leq x} a(n) e(\xi(n)).
\]
We are interested in the case where $a(n)$ is related to the normalized coefficients of a Hecke eigenform $f$ for the full modular group and $\xi(n)$ is a polynomial with real coefficients, although it is clear that other situations can be handled by the arguments in the paper. In particular, we consider the cases $a(n) = \lambda_f(n)$, the normalized Fourier coefficients of $f$, and $a(n)=\mu_f(n)$, the Dirichlet inverse of $\lambda_f(n)$. In both cases $a(n)$ is multiplicative and satisfies
\begin{equation}
\label{1-3}
|a(n)| \leq d(n),
\end{equation}
$d(n)$ being the divisor function.

\smallskip
There is a vast literature on estimates for $S_a(x,\xi)$, starting with the classical bounds for the linear case, where $\xi(n) = \alpha n$ with $\alpha\in\RR$; see e.g. Perelli \cite{Per/1984}, Jutila \cite{Jut/1987} and Fouvry-Ganguly \cite{Fo-Ga/2014}. In this paper we investigate some nonlinear cases.  When
\begin{equation}
\label{1-4}
\xi(n) = \sum_{\nu=0}^N a_\nu n^{\kappa_\nu}, \quad \kappa_0>\dots>\kappa_N>0, \quad a_\nu\in\RR
\end{equation}
and $\kappa_0\leq 1/2$, certain smoothed versions $\widetilde{S}_{\lambda_f}(x,\xi)$ of $S_{\lambda_f}(x,\xi)$ are well understood as special cases in the framework of the theory of nonlinear twists of $L$-functions developed by Kaczorowski-Perelli in a series of papers. Moreover, the same theory gives information on $\widetilde{S}_{\lambda_f}(x,\xi)$ for certain families of functions $\xi(n)$ of type \eqref{1-4} with leading exponent $\kappa_0>1/2$. We refer to Kaczorowski-Perelli \cite{Ka-Pe/resoI},\cite{Ka-Pe/resoII} for these results; see also Jutila \cite{Jut/1987a}. However, in the highly structured case where $\xi(n)$ is a polynomial of degree $k$, non-trivial bounds for $S_{\lambda_f}(x,\xi)$ or $\widetilde{S}_{\lambda_f}(x,\xi)$ are treated in the literature only when $k=2$; see Pitt \cite{Pit/2001} and few other papers stemming from it. Indeed, it is apparently difficult to proceed to higher degrees by the kind of arguments used in \cite{Pit/2001}, as these depend on delicate estimates involving sums of twisted Kloosterman sums. Moreover, at present general polynomials escape the analysis in \cite{Jut/1987a},\cite{Ka-Pe/resoI} and \cite{Ka-Pe/resoII}.

\smallskip
Although the bounds for $S_{\lambda_f}(x,\xi)$ in the nonlinear cases reported above show a power saving, it is nevertheless interesting to get weaker, but non-trivial, results for polynomials of arbitrary degree $\xi(n)$ and coefficients $\lambda_f(n)$ and $\mu_f(n)$.

\medskip
{\bf Theorem 2.} {\sl Let $P(n)$ be a polynomial with real coefficients and degree $k$. Then
\[
S_{\lambda_f}(x,P) \ll x\frac{\log\log x}{\log x} \quad \text{and} \quad S_{\mu_f}(x,P) \ll x\frac{\log\log x}{\sqrt{\log x}},
\]
where the constants in the $\ll$-symbols depend only on $f$ and $k$.}

\medskip
It will be clear from the proof that definitely better bounds can be obtained when the coefficients of $P(n)$ satisfy certain diophantine properties; see Section 3.2.

\medskip
In order to have the correct meaning of non-trivial bounds in the present case, we recall that 
\begin{equation}
\label{1-5}
\frac{x}{\log^\alpha x} \ll \sum_{n\leq x} |\lambda_f(n)| \ll \frac{x}{\log^\beta x}
\end{equation}
with $\alpha =0.211...$ and $\beta=0.118...$, see Wu \cite{Wu/2009}, while Elliott-Moreno-Shahidi \cite{E-M-S/1984} have shown that
\begin{equation}
\label{1-6}
\sum_{n\leq x} |\lambda_f(n)| \sim c\frac{x}{\log^\gamma x},
\end{equation}
with a certain constant $c=c(f)>0$ and $\gamma = 1-8/(3\pi) = 0.151...$, under the assumption of a strong form of the Sato-Tate conjecture. The referee pointed out that the known form of the Sato-Tate conjecture should imply at least that
\[
\sum_{n\leq x} |\lambda_f(n)| = \frac{x}{(\log x)^{\gamma+ o(1)}},
\]
since the distribution of $|\lambda_f(p)|$ is understood very well. Similar estimates hold for $|\mu_f(n)|$ as well.

\smallskip
Since the bounds in Theorem 2 are smaller than the left hand side of \eqref{1-5}, and hence than the right hand side of \eqref{1-6} as well, we may regard Theorem 2 as a quantitative form of orthogonality of $\lambda_f(n)$ and $\mu_f(n)$ to the exponentials $e(P(n))$. Moreover, Theorems 1 and 2 suggest the possibility of an extension of Sarnak's M\"obius Randomness conjecture  \cite{Sar/2011} to more general M\"obius functions, namely the Dirichlet coefficients of $1/L(s)$ for a suitable class of $L$-functions $L(s)$. A candidate for such a class are the primitive automorphic $L$-functions, of which the Hecke $L$-functions $L(s,f)$ are simple examples. For example, thanks to Theorem 1 some of the randomness results, already known for $\mu(n)$ via the BSZ theorem, should be transformable into randomness results for $\mu_f(n)$ in a rather direct way.

\smallskip
A major support to the M\"obius Randomness conjecture is provided by the fact that it follows from  the, a priori unrelated, Chowla conjecture; see \cite{Sar/2011}. One could therefore set up suitable extensions of these two conjectures and see if a similar implication holds between such extensions. However, this is apparently more tricky. Indeed, choosing for example $\mu_f(n)$ as a replacement of $\mu(n)$, a non-trivial bound for the extended M\"obius Randomness conjecture requires a saving of, roughly, $\log^\gamma x$ as in \eqref{1-6}. This adds some potential difficulties to be faced in such a procedure.

\medskip
{\bf Acknowledgements.} 
We wish to thank Sandro Bettin and Sary Drappeau for suggesting the use of the results by Shiu \cite{Shi/1980} and Nair \cite{Nai/1992} in the proof of Lemma 3.1. We also thank the referee for carefully reading our manuscript, and for pointing out several inaccuracies and improving the presentation at some points. This research was partially supported by the MIUR grant PRIN-2015 {\sl ``Number Theory and Arithmetic Geometry''}. The authors are members of the groups GNAMPA and GNSAGA of the Istituto Nazionale di Alta Matematica.

\bigskip
\section{Proof of Theorem 1}

\smallskip
We always assume that $x$ is sufficiently large.

\medskip
{\bf 2.1. Set up.}
For simplicity we assume that $H$ and $K$ in \eqref{1-1} are integers. Let $\nu\in[H,K]$,
\[
\I_\nu = ((\nu-1)^2,\nu^2], \quad \ P_\nu = \prod_{(H-1)^2 < p \leq \nu^2} p,
\]
\[
\P_\nu = \{p\in \I_\nu\}, \quad \M_\nu = \big\{m\in \big[1,\frac{x}{\nu^2}\big]: (m, P_\nu)=1 \big\},
\]
\[
 \P_\nu\M_\nu = \{pm: p\in\P_\nu, m\in\M_\nu\},
\]
\[
\I= \bigcup_{H\leq \nu \leq K} \P_\nu\M_\nu \ \ \text{and} \ \ \J= [1,x] \setminus \I;
\]
intervals are always meant as subsets of $\NN$. Note that each $n\in \P_\nu\M_\nu$ can be written in a unique way as $n=pm$ with $p\in\P_\nu$ and $m\in\M_\nu$, hence $|\P_\nu\M_\nu| = |\P_\nu||\M_\nu|$, and that $\P_\nu\M_\nu\subset [1,x]$. Moreover, the sets $\P_\nu\M_\nu$ are pairwise disjoint for $H\leq \nu\leq K$. 

\smallskip
Later on we will need certain bounds related to the sets $\P_\nu$, for $H\leq \nu\leq K$. Clearly, in view of the definition of $\P_\nu$, the Brun-Titchmarsh theorem immediately implies that
\begin{equation}
\label{2-1}
|\P_\nu| \ll \frac{\nu}{\log \nu}.
\end{equation}
Moreover, since by \eqref{1-1} we have $K^2 \leq x^{2\delta}$, a standard sieve estimate gives
\begin{equation}
\label{2-2}
\big|\big\{n\in[1,x]: n \ \text{has no prime factors in $\bigcup_{H\leq \nu \leq K}\P_\nu$}\big\} \big| \ll x \prod_{p|P_K}\big(1-\frac{1}{p}\big),
\end{equation}
see e.g. Theorem 3.5 of Halberstam-Richert \cite{Ha-Ri/1974}, and by Mertens' theorem we have
\begin{equation}
\label{2-3}
\prod_{p|P_\nu}\big(1-\frac{1}{p}\big) \ll \frac{\log H}{\log \nu}.
\end{equation}

\smallskip
Finally, we split $S(x)$ in \eqref{1-2} as
\begin{equation}
\label{2-4}
S(x) = \sum_{n\in\I} a(n) \phi(n) + \sum_{n\in\J} a(n) \phi(n) = S_\I(x) + S_\J(x),
\end{equation}
say.

\medskip
{\bf 2.2. Estimating $S_\I(x)$.}
We write
\begin{equation}
\label{2-5}
S_\I(x) = \sum_{H\leq \nu\leq K} \Big(\sum_{pm\in\P_\nu\M_\nu} a(pm)\phi(pm)\Big)  = \sum_{H\leq \nu\leq K} S_{\I,\nu},
\end{equation}
say. If $pm\in\P_\nu\M_\nu$ we have $(p,m)=1$, hence by the multiplicativity of $a(n)$, assumption (a) with the choice $y=x/\nu^2$ and $\overline{P}=P_\nu$, and \eqref{2-3} we get
\begin{equation}
\label{2-6}
\begin{split}
|S_{\I,\nu}| &= \Big|\sum_{m\in\M_\nu} a(m) \sum_{p\in\P_\nu} a(p)\phi(pm)\Big| \\
&\leq \Big(\sum_{m\in\M_\nu} |a(m)|^2\Big)^{1/2}  \Big(\sum_{m\in\M_\nu} \Big|\sum_{p\in\P_\nu} a(p)\phi(pm) \Big|^2 \Big)^{1/2} \\
&\ll  \Big(\frac{x\log H}{\nu^2\log \nu} \Big)^{1/2} \Big(\sum_{m\leq x/\nu^2} \Big|\sum_{p\in\P_\nu} a(p)\phi(pm) \Big|^2 \Big)^{1/2}.
\end{split}
\end{equation}
But thanks to assumption (b) with the choice $y=x/\nu^2$ and $\P=\P_\nu$, in view of \eqref{2-1}, $a(p)\ll1$ and $\phi(n) \ll 1$ we have
\begin{equation}
\label{2-7}
\begin{split}
\sum_{m\leq x/\nu^2} \Big|\sum_{p\in\P_\nu} a(p)\phi(pm) \Big|^2 &\ll \sum_{p,q\in\P_\nu} \Big| \sum_{m\leq x/\nu^2} \phi(pm)\overline{\phi(qm)} \Big| \\ 
&\ll  \frac{|\P_\nu| x}{\nu^2} + \sum_{\substack{p,q\in\P_\nu \\ p\neq q}} \Big| \sum_{m\leq x/\nu^2} \phi(pm)\overline{\phi(qm)} \Big| \\
&\ll \frac{x}{\nu\log \nu}\Big(1  + \frac{\tau\nu}{\log \nu} \Big),
\end{split}
\end{equation}
where $\tau = \tau(x)\leq 1$.

\smallskip
From \eqref{2-5},\eqref{2-6} and \eqref{2-7} we finally get
\begin{equation}
\label{2-8}
\begin{split}
S_\I(x) &\ll \sum_{H\leq \nu\leq K} \Big(\frac{x\log H}{\nu^2\log \nu} \Big)^{1/2} \Big(\frac{x}{\nu\log \nu} \big(1+\frac{\tau\nu}{\log \nu}\big)\Big)^{1/2} \\
&\ll x\sqrt{\log H} \left\{\sum_{H\leq \nu\leq K} \frac{1}{\nu^{3/2}\log \nu} + \sqrt{\tau} \sum_{H\leq \nu\leq K} \frac{1}{\nu \log^{3/2}\nu}\right\} \\
&\ll x \Big( \frac{1}{\sqrt{H\log H}} + \sqrt{\tau} \Big).
\end{split}
\end{equation}

\medskip
{\bf 2.3. Estimating $S_\J(x)$.}
We first define the following subsets of $[1,x]$:
\[
\begin{split}
\J_1^{(\nu)} &= \{n\in[1,x]: n \ \text{has exactly one prime divisor in $\P_\nu$ and none in $\bigcup_{H\leq h < \nu}\P_h$}\}, \\
\J_1 &= \bigcup_{H\leq \nu\leq K} \J_1^{(\nu)}, \\
\J_2 &= \big\{n\in[1,x]: n \ \text{has at least one prime factor in $\bigcup_{H\leq \nu \leq K}\P_\nu$}\big\},  \\
\J_3 &= \big\{n\in[1,x]: n \ \text{has no prime factors in $\bigcup_{H\leq \nu \leq K}\P_\nu$}\big\}.
\end{split}
\]
Clearly, $\J_1^{(\nu)} \supset \P_\nu\M_\nu$, hence $\J_1\supset \I$; moreover, $\J_2\cup \J_3 = [1,x]$ and $\J_2\cap \J_3 =\emptyset$. Thus, for future convenience, we write
\[
\J \subset (\J_1\setminus \I) \cup (\J_2\setminus \J_1) \cup \J_3.
\]
As a consequence, by assumption (a) with $P=1$ and $y=x$ we have that
\begin{equation}
\label{2-9}
\begin{split}
|S_\J(x)| &\ll \sum_{n\in\J_1\setminus\I} |a(n)| + \sum_{n\in\J_2\setminus\J_1} |a(n)| + \sum_{n\in\J_3} |a(n)| \\
&\ll x^{1/2}(|\J_1\setminus\I|^{1/2} + |\J_2\setminus\J_1|^{1/2}) + \Big(\sum_{n\in\J_3}|a(n)|^2\Big)^{1/2} |\J_3|^{1/2}.
\end{split}
\end{equation}

\smallskip
Clearly,
\[
\J_1^{(\nu)}\setminus \\P_\nu\M_\nu \subset \P_\nu \big(\frac{x}{\nu^2},\frac{x}{(\nu-1)^2}\big],
\]
hence by \eqref{2-1}
\begin{equation}
\label{2-10}
|\J_1\setminus\I| \ll \sum_{H\leq \nu\leq K} \frac{\nu}{\log \nu} \frac{x}{\nu^3} \ll \frac{x}{H\log H}.
\end{equation}
Moreover
\[
\J_2\setminus \J_1 \subset \bigcup_{H\leq \nu\leq K} \{n\in[1,x]: n \ \text{has at least two prime factors in $\P_\nu$}\},
\]
thus, again by \eqref{2-1},
\begin{equation}
\label{2-11}
|\J_2\setminus \J_1| \ll  \sum_{H\leq \nu\leq K} \sum_{p,q\in\P_\nu} \frac{x}{pq} \ll x  \sum_{H\leq \nu\leq K} \Big(\frac{|\P_\nu|}{(\nu-1)^2}\Big)^2 \ll \frac{x}{H\log^2 H}.
\end{equation}
Further, by assumption (a) with $y=x$ and $\overline{P}=P_K$, \eqref{2-2} and \eqref{2-3} we have
\begin{equation}
\label{2-12}
\sum_{n\in\J_3}|a(n)|^2 \ll x \frac{\log H}{\log K} \quad \text{and} \quad |\J_3| \ll x \frac{\log H}{\log K}.
\end{equation}

\smallskip
Collecting \eqref{2-9}-\eqref{2-12} we finally obtain that
\begin{equation}
\label{2-13}
S_\J(x) \ll x\Big(\frac{1}{\sqrt{H\log H}} +  \frac{\log H}{\log K}\Big),
\end{equation}
hence Theorem 1 follows from \eqref{2-4},\eqref{2-8} and \eqref{2-13}.

\bigskip
\section{Proof of Theorem 2}

\smallskip
We may clearly assume that the coefficients $\alpha_j$ of the polynomial $P(n)$ are reduced (mod 1). Hence, given large integers $Q_j=Q_j(x)>1$ for $1\leq j\leq k$, by Dirichlet's theorem  there exist $1\leq a_j \leq q_j\leq Q_j$ with $(a_j,q_j)=1$ such that
\begin{equation}
\label{3-1}
\Big|\alpha_j - \frac{a_j}{q_j}\Big| \leq \frac{1}{q_jQ_j}.
\end{equation}
Let $1< R_j <Q_j$, $R_j=R_j(x)$, be parameters to be chosen later on. With well established notation, we say that $\alpha_j$ belongs to the major arcs $\frak{M}_j$ if $\alpha_j$ satisfies \eqref{3-1} with some $1\leq q_j\leq R_j$, otherwise $\alpha_j$ belongs to the minor arcs $\frak{m}_j$. Moreover, with slight abuse of notation, we say that the polynomial $P(n)$ belongs to the major arcs $\frak{M}$ if $\alpha_j\in\frak{M}_j$ for every $j$, while $P(n)$ belongs to the minor arcs $\frak{m}$ if $\alpha_j\in\frak{m}_j$ for at least one $j$.

\smallskip
We treat these two cases for $P(n)$ by different techniques, but first we gather the required properties of the modular coefficients $\lambda_f(n)$ and $\mu_f(n)$, since the choice of the above parameters, as well as the quality of the final results, is heavily dependent on such properties.

\medskip
{\bf 3.1. Modular coefficients.}
We first list the results concerning $\lambda_f(n)$, starting with the well known bound given by the Ramanujan conjecture already recalled in \eqref{1-3}, namely
\begin{equation}
\label{3-2}
|\lambda_f(n)| \leq d(n).
\end{equation}
The next results are Theorem 1.3 of L\"u \cite{Lu/2009}, asserting that uniformly in $q$
\begin{equation}
\label{3-3}
\sum_{a=1}^q \Big| \sum_{\substack{n\leq x \\ n\equiv a \, (\text{mod $q$})}} \lambda_f(n) \Big|  \ll_f \sqrt{qx},
\end{equation}
and Jutila's theorem in \cite{Jut/1987}, according to which
\begin{equation}
\label{3-4}
\sum_{n\leq x} \lambda_f(n) e(\alpha n) \ll_f \sqrt{x}
\end{equation}
uniformly in $\alpha$. Moreover, it follows from the Rankin-Selberg convolution that
\begin{equation}
\label{3-5}
\sum_{n\leq x} |\lambda_f(n)|^2 \ll_f x,
\end{equation}
see Chapter 13 of Iwaniec \cite{Iwa/1997}. Let now $P= \prod_{z<p\leq w} p$.

\medskip
{\bf Lemma 3.1.} {\sl Let $P$ be as above with $z=z(x)\to\infty$ as $x\to\infty$ and $z < w < x$. Then}
\[
\sum_{\substack{n\leq x \\ (n,P)=1}} |\lambda_f(n)|^2 \ll_f x \prod_{p|P}\Big(1-\frac{1}{p}\Big).
\]

\medskip
{\it Proof.} Let $x$ be sufficiently large. Since $P$ depends on $x$, we consider the arithmetical function
\[
g_x(n) = 
\begin{cases}
|\lambda_f(n)|^2 & \text{if} \ (n,P)=1, \\
0 & \text{if} \ (n,P)>1.
\end{cases}
\]
Clearly, $g_x(n)$ is multiplicative and non-negative. Moreover, $g_x(n)$ belongs to the class $M=M(A_0,A_1)$, with certain $A_0,A_1$ independent of $x$, of multiplicative functions considered by Shiu \cite{Shi/1980} and Nair \cite{Nai/1992}; see p.259 of \cite{Nai/1992}. Indeed, from \eqref{3-2} we have $|g_x(p^\ell)| \leq d(p^\ell)^2 \leq (\ell+1)^2 \leq 4^\ell$ for every prime $p$ and $\ell\in\NN$, and \eqref{3-2} implies that there exists a function $c(\epsilon)>0$, independent of $x$, such that $g_x(n) \leq c(\epsilon) n^\epsilon$ for every $\epsilon>0$ and $n\in\NN$. Hence from the theorem on p.259 of \cite{Nai/1992} we get that
\begin{equation}
\label{3-6}
\sum_{\substack{n\leq x \\ (n,P)=1}} |\lambda_f(n)|^2 = \sum_{n\leq x} g_x(n) \ll x \prod_{p\leq x} \Big(1-\frac{1}{p}\Big) \exp\Big(\sum_{\substack{p\leq x \\ p\nmid P}} \frac{|\lambda_f(p)|^2}{p}\Big),
\end{equation}
the constant in the $\ll$-symbol being independent of $x$.

\smallskip
By \eqref{3-2} we have that
\[
\exp\Big(\sum_{\substack{p\leq x \\ p\nmid P}} \frac{|\lambda_f(p)|^2}{p}\Big) \asymp \prod_{p\leq x}  \Big(1+\frac{|\lambda_f(p)|^2}{p}\Big) \prod_{p|P}  \Big(1-\frac{|\lambda_f(p)|^2}{p}\Big).
\]
But the prime number theorem for $|\lambda_f(p)|^2$, see Rankin \cite{Ran/1973} or Perelli \cite{Per/1985} with $a=q=1$, implies that $|\lambda_f(p)|^2$ is asymptotically 1 on average, hence applying such a PNT three times, with $p\leq x$, $p\leq z$ and $p\leq w$, we finally obtain that
\begin{equation}
\label{3-7}
\exp\Big(\sum_{\substack{p\leq x \\ p\nmid P}} \frac{|\lambda_f(p)|^2}{p}\Big) \asymp  \prod_{p\leq x} \Big(1-\frac{1}{p}\Big)^{-1}  \prod_{p|P} \Big(1-\frac{1}{p}\Big).
\end{equation}
The lemma follows now from \eqref{3-6} and \eqref{3-7}. \qed

\smallskip
Now we turn to $\mu_f(n)$. We first note that from the Euler product for $L(s,f)^{-1}$ we have
\begin{equation}
\label{3-8}
\mu_f(n) = 
\begin{cases}
1 & \text{if} \ n=1 \\
(-1)^h\lambda_f(p_1\cdots p_h) & \text{if} \ n = p_1\cdots p_h(p_{h+1}\cdots p_r)^2, \ p_i\neq p_j \\
0 & \text{otherwise};
\end{cases}
\end{equation}
hence in particular from \eqref{3-2} we get
\begin{equation}
\label{3-9}
|\mu_f(p)| \leq 2.
\end{equation}
Next, the analogues of the bounds in \eqref{3-3} and \eqref{3-4} are given by the following lemmas.

\medskip
{\bf Lemma 3.2.} {\sl There exists an absolute constant $\delta_1>0$ such that, uniformly in $q$ and $1\leq a\leq q$, as $x\to\infty$ we have}
\[
 \sum_{\substack{n\leq x \\ n\equiv a \, (\text{mod $q$})}} \mu_f(n) \ll_f \sqrt{q} x e^{-\delta_1\sqrt{\log x}}.
\]

\medskip
{\it Proof.} The proof of this result is nowadays rather standard thanks to the non-existence of the Siegel zeros for the twisted Hecke $L$-functions associated with the cusp form $f$, proved by Hoffstein-Ramakrishnan \cite{Ho-Ra/1995} in 1995. Indeed, one may follow the arguments in Perelli \cite{Per/1984}, plugging in this extra information, or use those in Sections 4 and 7 of Fouvry-Ganguly \cite{Fo-Ga/2014}, already incorporating the Hoffstein-Ramakrishnan theorem. \qed

\medskip
{\bf Lemma 3.3.} {\sl There exists an absolute constant $\delta_2>0$ such that, uniformly in $\alpha$, as $x\to\infty$ we have}
\[
 \sum_{n\leq x} \mu_f(n) e(\alpha n) \ll_f  x e^{-\delta_2\sqrt{\log x}}.
\]

\medskip
{\it Proof.} Similarly as for the proof of Lemma 3.2. \qed

\smallskip
Finally, the analogues of \eqref{3-5} and Lemma 3.1 can be obtained as direct consequences by means of \eqref{3-8}. Indeed, for $\overline{P}=1$ or $\overline{P}=P$ as in Lemma 3.1 with $w\leq 2x^{2\delta}$, $\delta$ being as in \eqref{1-1}, from \eqref{3-5} and Lemma 3.1 we have
\begin{equation}
\label{3-10}
\begin{split}
\sum_{\substack{n\leq x \\ (n,\overline{P})=1}} |\mu_f(n)|^2 &= \sum_{\substack{p_1\cdots p_h (p_{h+1}\cdots p_r)^2 \leq x \\ p_j\nmid \overline{P}}} |\mu_f\big(p_1\cdots p_h (p_{h+1}\cdots p_r)^2\big)|^2 \\
&\leq \sum_{d\leq \sqrt{x}} \, \sum_{\substack{p_1\cdots p_h \leq x/d^2 \\ p_j\nmid \overline{P}}} |\lambda_f(p_1\cdots p_h)|^2 \\
&\leq  \sum_{d \leq \sqrt{x}}  \sum_{\substack{m\leq x/d^2 \\ (m,\overline{P})=1}} |\lambda_f(m)|^2 \\
&\ll_f  \sum_{d\leq  \frac{1}{\sqrt{2}} x^{(1-2\delta)/2}}  \sum_{\substack{m\leq x/d^2 \\ (m,\overline{P})=1}} |\lambda_f(m)|^2 + x^{(1+2\delta)/2} \\
&\ll_f x\prod_{p|\overline{P}} \Big(1-\frac{1}{p}\Big).
\end{split}
\end{equation}

\medskip
{\bf 3.2. Major arcs estimates.}
Recalling the notation after \eqref{3-1}, we start with the case where $P(n)$ belongs to $\frak{M}$. Clearly, the size of the $R_j$ will depend on the level of distribution of the coefficients $\lambda_f(n)$ and $\mu_f(n)$ in arithmetic progressions. We indeed have that
\[
P(n) = \sum_{j=1}^k \frac{a_j}{q_j} n^j + \sum_{j=1}^k \big(\alpha_j-\frac{a_j}{q_j} \big) n^j = \overline{P}(n) + R(n),
\]
say, and hence, denoting by $a(n)$ either $\lambda_f(n)$ or $\mu_f(n)$, by partial summation we get
\begin{equation}
\label{3-11}
\begin{split}
S_a(x,P) &:= \sum_{n\leq x} a(n) e(P(n)) = \sum_{n\leq x} a(n) e(\overline{P}(n) + R(n)) \\
&\ll |S_a(x,\overline{P})| + x \max_{1\leq j\leq k} \max_{1\leq t\leq x} \frac{t^{j-1}}{q_jQ_j}|S_a(t,\overline{P})|.
\end{split}
\end{equation}
Moreover, writing 
\[
\text{$q=$ lcm$ \, (q_1,\dots,q_k)$ and $\overline{P}(n)=\frac{1}{q} \sum_{j=1}^k b_j n^j :=\frac{1}{q} \widetilde{P}(n)$, $b_j\in\NN$}, 
\]
with obvious notation we obtain that
\begin{equation}
\label{3-12}
|S_a(t,\overline{P})| =\Big| \sum_{b=1}^q e(\widetilde{P}(b)/q) \Big(\sum_{\substack{n\leq t \\ n\equiv b \, (\text{mod $q$})}} a(n)\Big)\Big|  \leq \sum_{b=1}^q|S_a(t;q,b)|.
\end{equation}

\smallskip
{\bf Case 1}: $a(n) = \lambda_f(n)$. By \eqref{3-11},\eqref{3-12} and \eqref{3-3} we have
\begin{equation}
\label{3-13}
S_{\lambda_f}(x,P) \ll (qx)^{1/2} (1 + \max_{1\leq j\leq k} x^jQ_j^{-1}).
\end{equation}
In this case we choose
\begin{equation}
\label{3-14}
Q_j= x^{j-c_j} \quad \text{and} \quad R_j = x^{c'_j}
\end{equation}
with $c_1,\dots,c_k,c'_1,\dots,c'_k>0$, $c_j<1$ and $c'_j<j-c_j$ to be determined later on. Therefore, from the definition of $q$, \eqref{3-13} and \eqref{3-14}, if $P(n)$ belongs to $\frak{M}$ we obtain
\begin{equation}
\label{3-15}
S_{\lambda_f}(x,P) \ll x^{\gamma_1} \quad \text{with} \quad \gamma_1 = \frac12 + \max_{1\leq j\leq k} c_j +\frac12\sum_{j=1}^k c'_j.
\end{equation}

\smallskip
{\bf Case 2}: $a(n) = \mu_f(n)$. In this case we choose
\begin{equation}
\label{3-16}
Q_j = x^j e^{-\beta_j\sqrt{\log x}} \quad \text{and} \quad R_j =e^{\beta'_j\sqrt{\log x}},
\end{equation}
with $\beta_1,\dots,\beta_k,\beta'_1,\dots,\beta'_k>0$ to be determined later on. Thus from \eqref{3-11},\eqref{3-12}, Lemma 3.2, \eqref{3-16} and the definition of $q$, if $P(n)$ belongs to $\frak{M}$ we obtain
\begin{equation}
\label{3-17}
S_{\mu_f}(x,P) \ll x e^{-\gamma_1'\sqrt{\log x}} \quad \text{with} \quad  \gamma_1'= \delta_1-\max_{1\leq j\leq k} \beta_j - \frac{3}{2}\sum_{j=1}^k \beta'_j.
\end{equation}

\medskip
{\bf 3.3. A Weyl type lemma.}
In order to verify assumption (b) in Theorem 1 with our choice $\phi(n)=e(P(n))$, when $a(n)=\mu_f(n)$ we need a sharper version of the classical Weyl lemma on the bound for exponential sums with polynomial values; see Theorem 2 in Chapter 3 of Montgomery \cite{Mon/1994}. Essentially, we need to replace the term $x^\epsilon$ in the classical bound by a power of $\log x$, plus other minor variants. Actually, the result we need is in the spirit of the lemma on p.199 of Perelli-Zaccagnini \cite{Pe-Za/1995}; since we could not trace the required result in the literature, we provide a proof here. 

\smallskip
We first state a slight variant of a classical auxiliary lemma, whose proof follows closely that of (9) in Chapter 3 of \cite{Mon/1994}.

\medskip
{\bf Lemma 3.4.} {\sl Let  $|\alpha-a/q|\leq C/q^2$ with some $1\leq a<q$, $(a,q)=1$ and $C\geq1$, and let $M,N\geq1$. Then, writing $\|\xi\|$ for the distance of $\xi$ from the nearest integer, we have}
\[
\sum_{n=1}^N \min\Big(M,\frac{1}{\|\alpha n\|}\Big) \ll C \Big( \frac{MN}{q} + N\log q + M + q\log q\Big).
\]

\medskip
The next result gives the required form of Weyl's lemma.

\medskip
{\bf Lemma 3.5.} {\sl Let $d\geq 2$, $U(n) = \alpha n^d + \alpha_{d-1}n^{d-1}+ \cdots + \alpha_1n$ with $\alpha_j\in\RR$ and $\alpha$ as in Lemma $3.4$. Then, writing $\kappa = 2^{1-d}$, for any $Z > 1$ we have 
\[
W(y,U) := \sum_{n\leq y} e(U(n)) \ll y\Big(\frac{CZ}{q} + \frac{CZ}{y}\log q +CZ\frac{q\log q}{y^d} + \frac{\log^A y}{Z}   \Big)^\kappa,
\]
where $A=A(d)$ is  a certain constant and the constant in the $\ll$-symbol depends only on $d$.}

\medskip
{\it Proof.} We may suppose that $y\in\NN$; moreover, here we denote by $\tau_{\ell}(n)$ the $\ell$th divisor function. Following the proof of the above mentioned Theorem 2 in \cite{Mon/1994}, by Weyl's differencing method applied $d-1$ times we get
\begin{equation}
\label{3-18}
|W(y,U)|^{2^{d-1}} \ll y^{2^{d-1}-1} + y^{2^{d-1}-d} \sum_{h_1,\dots,h_{d-1}} \min\Big(y, \frac{1}{\|d!h_1\cdots h_{d-1} \alpha\|}\Big),
\end{equation}
where $h_j\in[1,y-1-h_{j-1}]$ (here $h_0=0$) and hence $d!h_1\cdots h_{d-1} \leq d!y^{d-1}$. Therefore we have that
\begin{equation}
\label{3-19}
 \sum_{h_1,\dots,h_{d-1}} \min\Big(y, \frac{1}{\|d!h_1\cdots h_{d-1} \alpha\|}\Big) \leq \sum_{h\leq d!y^{d-1}} \tau_{d-1}(h) \min\Big(y,\frac{1}{\|h\alpha\|}\Big).
\end{equation}

\smallskip
Let now $Z> 1$ and $\H_Z^-$ be the set of the $h\leq d!y^{d-1}$ with $\tau_{d-1}(h) \leq Z$, and $\H_Z^+=[1,d!y^{d-1}]\setminus \H_Z^-$. Thus from Lemma 3.4 we get
\begin{equation}
\label{3-20}
\sum_{h\in\H_Z^-} \tau_{d-1}(h) \min\Big(y,\frac{1}{\|h\alpha\|}\Big) \ll CZ \Big( \frac{y^d}{q} + y^{d-1}\log q + y + q\log q\Big),
\end{equation}
while recalling the standard bounds for the mean-square of the $(d-1)$th divisor function we obtain
\begin{equation}
\label{3-21}
\begin{split}
\sum_{h\in\H_Z^+} \tau_{d-1}(h) \min\Big(y,\frac{1}{\|h\alpha\|}\Big) &\ll \frac{1}{Z} \sum_{h\leq d!y^{d-1}} \tau_{d-1}(h)^2 \min\Big(y,\frac{1}{\|h\alpha\|}\Big) \\
&\ll \frac{y}{Z}  \sum_{h\leq d!y^{d-1}} \tau_{d-1}(h)^2 \ll \frac{y^d}{Z} \log^cy
\end{split}
\end{equation}
with a certain $c=c(d)$. The result follows now from \eqref{3-18}-\eqref{3-21}, since $d\geq 2$. \qed

\smallskip
We finally recall that, under the same hypotheses of Lemma 3.5, the standard Weyl bound becomes
\begin{equation}
\label{3-22}
W(y,U) \ll y^{1+\epsilon} C^\kappa \Big(\frac{1}{q} + \frac{1}{y} +\frac{q}{y^d} \Big)^\kappa \quad \text{for every} \ \epsilon>0.
\end{equation}

\medskip
{\bf 3.4. Minor arcs estimates.}
Finally, again recalling the notation after \eqref{3-1}, we deal with the case where $P(n)$ belongs to $\frak{m}$. In this case our basic tool will be Theorem 1, with the choice of $a(n)$ as in Section 3.2, i.e. either $\lambda_f(n)$ or $\mu_f(n)$, and $\phi(n)=e(P(n))$. Thus we have to show that the assumptions in Theorem 1 are satisfied with such choices. Again we consider separately the two cases of $a(n)$, but first we proceed to some preliminary reductions common to both cases. Let
\[
d= \max \,\{1\leq j\leq k: q_j>R_j\}.
\]

\smallskip
Suppose first that $d=1$; in this case we argue directly, without appealing to Theorem 1 nor to Lemma 3.5. Recalling \eqref{3-1}, \eqref{3-11} and that an empty sum equals 0, we write
\[
P(n) =  \alpha_1 n + \sum_{j=2}^k \frac{a_j}{q_j} n^j + \sum_{j=2}^k \Big(\alpha_j-\frac{a_j}{q_j}\Big) n^j = L(n) + R_1(n) + R_2(n),
\]
say, hence arguing as in Section 3.2, by partial summation we get
\begin{equation}
\label{3-23}
S_a(x,P) \ll |S_a(x,L+R_1)| + x \max_{2\leq j\leq k} \max_{1\leq t\leq x} \frac{t^{j-1}}{q_jQ_j}  |S_a(t,L+R_1)|.
\end{equation}
Moreover, writing
\[
\overline{q}= \ \text{lcm} \, (q_2,\dots,q_k) \quad \text{and} \quad R_1(n)=\frac{1}{\overline{q}} \sum_{j=2}^k A_j n^j :=\frac{1}{\overline{q}} \widetilde{R}_1(n), \ A_j\in\NN, 
\]
thanks to the orthogonality of additive characters we have
\begin{equation}
\label{3-24}
\begin{split}
S_a(t,L+R_1) &=  \sum_{b=1}^{\overline{q}} e(\widetilde{R}_1(b)/\overline{q}) \Big(\sum_{\substack{n\leq t \\ n\equiv b \, (\text{mod $\overline{q}$})}} a(n) e(L(n)) \Big) \\
&= \sum_{b=1}^{\overline{q}} e(\widetilde{R}_1(b)/\overline{q}) \frac{1}{\overline{q}} \sum_{c=1}^{\overline{q}} e(-bc/\overline{q}) \sum_{n\leq t} a(n) e\big((\alpha_1+c/\overline{q})n\big) \\
&\ll \overline{q} \max_{\alpha\in[0,1]} \Big| \sum_{n\leq t} a(n) e(\alpha n) \Big|.
\end{split}
\end{equation}

\smallskip
Suppose now that $2\leq d\leq k$; in this case we use both Theorem 1 and Lemma 3.5. Given $\P$ as in assumption (b) of Theorem 1 and $p,q\in\P$ with $p\neq q$, writing $C_j=p^j-q^j \ll z^j$ and recalling that $\phi(n)=e(P(n))$ we have that
\[
\phi(pm)\overline{\phi(qm)} = e\Big(\sum_{j=1}^k C_j\alpha_jm^j\Big).
\]
Arguing similarly as before we split the above polynomial as
\[
\begin{split}
\sum_{j=1}^k C_j\alpha_jm^j &=  \sum_{j=1}^{d} C_j\alpha_jm^j +\sum_{j=d+1}^k C_j\frac{a_j}{q_j}m^j + \sum_{j=d+1}^k C_j\Big(\alpha_j-\frac{a_j}{q_j}\Big)m^j \\
&= U(m) + V(m) + \widetilde{R}(m),
\end{split}
\]
say. Thus, writing
\[
W(y,U+V) = \sum_{m\leq y} e(U(m)+V(m)),
\]
by partial summation we get
\begin{equation}
\label{3-25}
\sum_{m\leq y} \phi(pm)\overline{\phi(qm)}  \ll |W(y,U+V)| + y \max_{d+1\leq j\leq k} \max_{1\leq t\leq y} \frac{t^{j-1}}{q_jQ_j}  |W(t,U+V)|.
\end{equation}
Moreover, letting this time $\widetilde{q} = \text{lcm} \, (q_{d+1},\dots,q_k)$, arguing as for \eqref{3-24} we obtain
\begin{equation}
\label{3-26}
W(t,U+V) \ll \widetilde{q} \max_{b=1,\dots,\widetilde{q}}\Big| \sum_{n\leq t} e(U(n) + (b/\widetilde{q})n)\Big|.
\end{equation}
But, since $U(n) + (b/\widetilde{q})n$ has degree $d\geq2$, we may apply Lemma 3.5 or \eqref{3-22} to the right hand side of \eqref{3-25}. Hence in view of \eqref{3-22} with $y=t$, $\alpha=\alpha_d$, $C=C_d \ll z^d$ and $q=q_d$ with $R_d<q_d\leq Q_d$, from the definition of $d$ and $\widetilde{q}$, \eqref{3-25} and \eqref{3-26} we get, after taking the maximum over $1\leq t\leq y$, that
\begin{equation}
\label{3-27}
\sum_{m\leq y} \phi(pm)\overline{\phi(qm)} \ll y^{1+\epsilon} z^{\kappa d} R_{d+1} \cdots R_k \Big(1+\max_{d+1\leq j\leq k} \frac{y^j}{Q_j}\Big) \Big(\frac{1}{R_d} + \frac{1}{y} + \frac{Q_d}{y^d}\Big)^\kappa
\end{equation}
with $\kappa = 2^{1-d}$. Alternatively, appealing instead to Lemma 3.5 with the same choices as above, again from the definition of $d$ and $\widetilde{q}$, \eqref{3-25} and \eqref{3-26}, arguing as before we have
\begin{equation}
\label{3-28}
\begin{split}
\sum_{m\leq y} \phi(pm)\overline{\phi(qm)} &\ll y R_{d+1} \cdots R_k \Big(1+\max_{d+1\leq j\leq k} \frac{y^j}{Q_j}\Big) \\
&\times  \Big(\frac{z^dZ}{R_d} + \frac{z^dZ}{y}\log Q_d + z^dZ\frac{Q_d\log Q_d}{y^d} + \frac{\log^A y}{Z}\Big)^\kappa,
\end{split}
\end{equation}
with any $Z>1$ and still $\kappa = 2^{1-d}$.

\smallskip
{\bf Case 1}: $a(n) = \lambda_f(n)$. We first deal with the case $d=1$. From \eqref{3-4}, the definition of $\overline{q}$, \eqref{3-14},\eqref{3-23} and \eqref{3-24}, for $d=1$ and $P\in\frak{m}$ we get
\begin{equation}
\label{3-29}
S_{\lambda_f}(x,P) \ll \overline{q} x^{1/2} (1+\max_{2\leq j\leq k} x^jQ_j^{-1}) \ll x^{\gamma_2} \quad \text{with} \quad \gamma_2 = \frac12+\max_{2\leq j\leq k} c_j + \sum_{j=2}^k c'_j.
\end{equation}

\smallskip
For $d\geq 2$ we use Theorem 1, thus we have to verify its assumptions. Clearly $\lambda_f(p)\ll1$ follows from \eqref{3-2}, while assumption (a) follows from \eqref{3-5} and Lemma 3.1, without imposing any condition on $H$ and $K$ in addition to \eqref{1-1}. Concerning assumption (b), from \eqref{3-27} and \eqref{3-14} we have that
\begin{equation}
\label{3-30}
\sum_{m\leq y} \phi(pm)\overline{\phi(qm)} \ll \tau y \quad \text{with} \quad y\asymp x/z
\end{equation}
is satisfied uniformly for $p,q$ as in (b), $p\neq q$, with the choice
\begin{equation}
\label{3-31}
\tau = x^\epsilon z^{\kappa d} x^{c'_{d+1}+\cdots +c'_k} \Big(\frac{1}{x^{\kappa c'_d}} + \big(\frac{z}{x})^\kappa + (\frac{z^d}{x^{c_d}}\big)^\kappa \Big) \big(1+ \max_{d+1\leq j\leq k} z^{-j}x^{c_j}\big).
\end{equation}
Hence, choosing $\delta$ in \eqref{1-1} sufficiently small, since $z\ll x^{2\delta}$ we have that \eqref{3-30} holds with
\begin{equation}
\label{3-32}
\tau= x^{-c_0}
\end{equation}
with a small constant $c_0>0$, depending on $\epsilon$, $\delta$ and the various constants involved in \eqref{3-31}, provided
\begin{equation}
\label{3-33}
c'_{d+1}+\cdots +c'_k + \max_{d+1\leq j\leq k} c_j < \min (\kappa,\kappa c_d,\kappa c'_d).
\end{equation}
In order to avoid a simple but tedious optimization, we now observe that clearly \eqref{3-33} holds if all constants $c_j$ and $c'_j$ are chosen sufficiently small and satisfying, for example,
\[
c_{j+1} \leq 2^{-10j} c_j \quad \text{and} \quad c'_{j+1} \leq 2^{-10j} c'_j \quad \text{for} \quad 1\leq j\leq k-1.
\]

\smallskip
Therefore, after a trivial summation over $p\neq q$, we have that assumption (b) is satisfied with the choice of $\tau$ in \eqref{3-32}, again without imposing any condition on $H$ and $K$ in addition to \eqref{1-1}. Thus from Theorem 1 we obtain that
\[
S_{\lambda_f}(x,P) \ll x\Big(\frac{1}{\sqrt{H\log H}} + x^{-c_0} + \frac{\log H}{\log K}\Big),
\]
hence, choosing for example $H= \log^2x$ and $K = x^\delta$, for $d\geq 2$ we get
\begin{equation}
\label{3-34}
S_{\lambda_f}(x,P) \ll x\frac{\log\log x}{\log x}.
\end{equation}
Finally, since with the above choice of the constants $c_j$ and $c'_j$ we also have that the constants $\gamma_1$ and $\gamma_2$ in \eqref{3-15} and \eqref{3-29} are both $<1$, the first assertion of Theorem 2 follows from \eqref{3-15}, \eqref{3-29} and \eqref{3-34}.

\smallskip
{\bf Case 2}: $a(n) = \mu_f(n)$. The deduction of the second assertion of Theorem 2 is similar, so we give only a brief account of the needed changes. From Lemma 3.3, the definition of $\overline{q}$, \eqref{3-16},\eqref{3-23} and \eqref{3-24}, for $d=1$ and $P$ belongs to $\frak{m}$ we get
\begin{equation}
\label{3-35}
S_{\mu_f}(x,P)  \ll xe^{-\gamma'_2\sqrt{\log x}} \quad \text{with} \quad \gamma'_2 = \delta_2 - \max_{2\leq j\leq k} \beta_j - \sum_{j=2}^k \beta'_j.
\end{equation}

\smallskip
For $d\geq 2$ we use again Theorem 1. Also in this case, thanks to \eqref{3-9} and \eqref{3-10}, $\mu_f(p)\ll 1$ and assumption (a) are satisfied without imposing any condition on $H$ and $K$ in addition to \eqref{1-1}. In order to verify assumption (b), this time we use \eqref{3-28} and \eqref{3-16} to obtain that \eqref{3-30} is satisfied uniformly for $p,q$ as in (b), $p\neq q$, with the choice (here we write $L=\sqrt{\log x}$)
\begin{equation}
\label{3-36}
\begin{split}
\tau = (z^{d}Z)^{\kappa} &e^{(\beta'_{d+1}+\cdots +\beta'_k)L} \Big(e^{-\kappa \beta'_dL} + \big(\frac{z\log x}{x}\big)^\kappa  + \big(z^{d} e^{-\beta_dL}\log x\big)^\kappa \Big) \big(1+\max_{d+1\leq j\leq k} z^{-j} e^{\beta_jL}\big) \\
&+ e^{(\beta'_{d+1}+\cdots +\beta'_k)L} \Big( \frac{\log^{A}x}{Z}\Big)^\kappa \big( 1 + \max_{d+1\leq j\leq k} z^{-j} e^{\beta_jL} \big).
\end{split}
\end{equation}
Assuming that
\begin{equation}
\label{3-37}
K = e^{\delta\sqrt{\log x}} \quad \text{and} \quad Z= e^{\mu\sqrt{\log x}},
\end{equation}
and hence $z \leq 2e^{2\delta\sqrt{\log x}}$, we see that the dependence on the constants $\beta_j$, $\beta'_j$, $\delta$ and $\mu$ in \eqref{3-36} is structurally very similar to that in \eqref{3-31}. Hence similar arguments as before show that there exists a choice of the involved constants such that \eqref{3-30} holds with the choice
\begin{equation}
\label{3-38}
\tau = e^{-c'_0\sqrt{\log x}},
\end{equation}
where $c'_0>0$ is a small constant. Therefore, in view of \eqref{3-37} and \eqref{3-38}, choosing for example $H= \log x$ and $K = e^{\delta\sqrt{\log x}}$ in Theorem 1, for $d\geq 2$ we get
\begin{equation}
\label{3-39}
S_{\mu_f}(x,P) \ll x\frac{\log\log x}{\sqrt{\log x}}.
\end{equation}
Moreover, with such choices of the constants we also have that the values of $\gamma'_1$ and $\gamma'_2$ in \eqref{3-17} and \eqref{3-35} are both $>0$, and the second assertion of Theorem 2 follows from \eqref{3-17}, \eqref{3-35} and \eqref{3-39}. The proof is now complete.

\bigskip

\ifx\undefined\bysame{poly}.
\newcommand{\bysame}{\leavevmode\hbox to3em{\hrulefill}\ ,}
\fi

\bigskip
\bigskip
\bigskip
\noindent

\medskip
\noindent
Mattia Cafferata, Dipartimento di Scienze Matematiche, Fisiche e Informatiche, Universit\`a di Parma, Parco Area delle Scienze 53/a, 43124 Parma, Italy;  \url{mattia.cafferata@unife.it}

\medskip
\noindent
Alberto Perelli, Dipartimento di Matematica, Universit\`a di Genova, via Dodecaneso 35, 16146 Genova, Italy; \url{perelli@dima.unige.it}

\medskip
\noindent
Alessandro Zaccagnini, Dipartimento di Scienze Matematiche, Fisiche e Informatiche, Universit\`a di Parma, Parco Area delle Scienze 53/a, 43124 Parma, Italy; 

\noindent
\url{alessandro.zaccagnini@unipr.it}

\end{document}